

%


\documentstyle{amsppt}
\magnification=1200
\NoBlackBoxes
\loadeusm

\topmatter
\title UMBRELLAS AND POLYTOPAL APPROXIMATION
OF THE EUCLIDEAN BALL
\endtitle
\author
Yehoram Gordon,
Shlomo Reisner, and
Carsten Sch\"utt
\endauthor
\address
Yehoram Gordon,
Department of Mathematics, Technion, 32000 Haifa, Israel
\endaddress
\thanks Yehoram Gordon was partially supported by the Fund for the 
Promotion of Research at the Technion
\endthanks
\address
Shlomo Reisner, Department of Mathematics and School of
Education--Oranim, University of Haifa, Israel,   
and
Department of Mathematics and Computer Science,
University of Denver, Denver, Colorado 80208, USA
\endaddress
\address 
Carsten Sch\"utt,
Oklahoma State University,
Department of Mathematics,
Stillwater, Oklahoma 74078, USA, and
Christian-Albrechts Universit\"at,
Mathematisches Seminar, 24098 Kiel, Germany
\endaddress
\thanks Carsten Sch\"utt has been supported by NSF-grant DMS-9301506
\endthanks
\subjclass 52A22
\endsubjclass
\abstract
There are two positive, absolute constants $c_{1}$ and $c_{2}$ so that
the volume of the difference set of the $d$-dimensional Euclidean
ball and an inscribed polytope with n vertices is larger than
$$
c_{2}\ d\ {n}^{-\frac{2}{d-1}}vol_d(B^d_2)
$$
for $n \geq (c_{1}\ d)^{\frac{d-1}{2}}$. 

\endabstract
\endtopmatter 
\document
\break

We study here the approximation of a convex body in $\Bbb R^{d}$
by a polytope with at most n vertices. There are many means to measure
the approximation, the two most common are the
Hausdorff distance or the symmetric difference metric. The Hausdorff
distance between two convex bodies $K$ and $C$ is
$$
d_{H}(K,C)=\max\{\max_{x \in C} \min_{y \in K} \|x-y\|_{2},
\max_{y \in K} \min_{x \in C} \|x-y\|_{2} \}
$$
where $\|x\|_{2}$ is the Euclidean norm of $x$.
The symmetric difference metric is the volume of the difference set.
$$
d_{S}(K,C)=vol_{d}(K \triangle C).
$$
Bronshtein and Ivanov [BI] and Dudley [$\text{D}_{1}, \text{D}_{2}$]
showed that for every convex body there is a constant $c$ such that
for every $n$ there is a polytope $P_{n}$ with at most $n$ vertices and
$$
d_{H}(K,P_{n}) \leq c n^{-\frac{2}{d-1}}.
$$
This can be used to show the same estimate for the symmetric
difference metric. Gruber and Kenderov [GK] showed that the inverse
inequality holds if $K$ has a $C^{2}$-boundary:
$$
d_{S}(K,P_{n}) \geq c n^{-\frac{2}{d-1}}.
$$
Macbeath [Mac] showed that the approximation of a convex body is always
better than that of the Euclidean sphere.
Gruber [$\text{Gr}_{2}$] obtained an asymptotic formula. If a convex body 
$K$ in $\Bbb{R}^{d}$ has a $C^{2}$-boundary 
with everywhere positive curvature, then we have
$$
\inf \{d_{S}(K,P_{n})|P_{n} \subset K\ \text{and $P_{n}$ \ has at most
n vertices} \}
\sim
$$
$$
\frac{1}{2} \text{del}_{d-1} \int_{\partial K} \kappa(x)^{\frac{1}{d+1}}d\mu(x) 
(\frac{1}{n})^{\frac{2}{d-1}}
$$
where $del_{d-1}$ is a constant that is connected with triangulations.
In [GMR$_{1}$, GMR$_{2}$] it was shown constructively that for all dimensions
$d$, all convex
bodies $K$, and all $n \geq 2$ there is a polytope $P_{n}$ with $n$ vertices 
that is contained in $K$ such that
$$
vol_{d}(K)-vol_{d}(P_{n}) \leq c\ d\ vol_{d}( K) n^{-\frac{2}{d-1}}
$$
where $c$ is a numerical constant. This estimate can also be derived from
[BI] and [D$_1$, D$_2$]. So the question was whether the factor $d$
was necessary, or, in other words, what is the order of magnitude of the
constant
$del_{d}$. The result in this paper shows that there are absolute constants 
$c_{1}$ and $c_{2}$ with
$$
c_{1} \leq \text{del}_{d} \leq  c_{2}.
$$ 
In fact, we have
$$
\text{del}_{d-1} \leq \frac{32}{7}(\frac{vol_{d-1}(\partial B_{2}^{d})}
{vol_{d-1}(B_{2}^{d-1})})^{\frac{2}{d-1}}.
$$
This follows from estimate (1) below.
\par 
In this paper we want to show that the volume of the difference 
set of the $d$-dimensional Euclidean
ball and an inscribed polytope with n vertices is 
larger than
$$
cd\ vol_{d}(B_{2}^{d})\ {n}^{-\frac{2}{d-1}}
$$
We want to reduce the computation of the volume of the difference set
to that of the following set: The set between a $d-1$ dimensional face
of the polytope and the boundary of the sphere.
Intuitively it is clear that the faces should be simplices and that the polytope
should have rather regular features. This leads us to the assumption
that the volume of the set between a $d-1$ dimensional face
of the polytope and the boundary of the sphere equals in average approximately
the surface area of the face times the height of the cap of the Euclidean ball 
that is determined by that face.
\par 
There are two technical difficulties. The number
of faces does not necessarily correspond to the number of vertices.
In fact, a heuristic argument shows that the number of faces is
of the order of the number of vertices times $d^{\frac{d}{2}}$.
Secondly, although we may assume that the faces are simplices, we
may not assume that they are regular or close to regular. This
is expressed in the following way. If $F$ is a face and $H$ the 
hyperplane containing $F$ then the distance of the centers of
gravity of $F$ and $H \cap B_{2}^{d}$ may be large.
\par  
Hyperplanes are usually denoted by $H$ and the closed
halfspaces associated with $H$ by $H^{+}$ and $H^{-}$.
$H(x,\xi)$ is the hyperplane that passes through $x$
and is orthogonal to $\xi$.
\par
The $d-1$ dimensional faces of a polytope in $\Bbb R^{d}$
are denoted by $F_{j}$. The hyperplanes containing $F_{j}$
are denoted by $H_{j}$. $H_{j}^{+}$ denotes the halfspace
containing $P$. 
\par
For a polytope $P$ that is contained
in $B_{2}^{d}$ the height or width of $B_{2}^{d} \cap
H_{j}^{-}$ is $h_{j}$ and the radius of $B_{2}^{d} \cap
H_{j}$ is $r_{j}$.
\par
$cg(M)$ is the center of gravity of the set $M$.
\par
$[A,B]$ denotes the convex hull of the sets $A$ and $B$.
The radial projection $rp(M)$ of a set $M$ in $B_{2}^{d}$ is
$$
rp(M)=\{ \xi \in \partial B_{2}^{d} | [0,\xi] \cap M \neq \emptyset \}.
$$
\vskip 10mm

\proclaim{\smc Theorem 1}
There are two positive constants  $c_{1}$ and $c_{2}$ so that we
have for all $d,d \geq 2$, and all $n,n \geq (c_{1}\ d)^{\frac{d-1}{2}}$, and
all polytopes
$P_{n}$ that are contained in the Euclidean unit ball $B_{2}^{d}$ and have 
$n$ vertices
$$
vol_{d}(B_{2}^{d})-vol_{d}(P_{n}) \geq c_{2} \ d\
vol_{d}(B_{2}^{d}) n^{-\frac{2}{d-1}}.
$$ 
\endproclaim
\vskip 3mm

In particular we have by Theorem 1 that there are positive constants 
$c_3$ and $c_4$ such that 
$$
vol_{d}(B_{2}^{d})-vol_{d}(P_{n}) \geq c_4 vol_{d}(B_{2}^{d}) 
$$
if $n \leq (c_3 d)^{\frac{d-1}{2}}$.

\proclaim{\smc Lemma 2}
(i) For all $x,0<x,$ there is a $\theta,0< \theta <1,$ such that
$$
\Gamma(x+1)=\sqrt{2\pi}x^{x+\frac{1}{2}}
\exp(-x+\frac{\theta}{12x}).
$$
(ii)
$$
vol_{d}(B_{2}^{d})= \frac{\pi^{\frac{d}{2}}}{\Gamma(\frac{d}{2}+1)} 
\leq 
\frac{\pi^{\frac{d-1}{2}}(2e)^\frac{d}{2}}{d^{\frac{d+1}{2}}}.
$$
\endproclaim
\vskip 3mm

The following lemma is due to Bronshtein and Ivanov [BI] and
Dudley $[D_{1}, D_{2}]$.

\proclaim{\smc Lemma 3}
For all dimensions $d$, $d \geq 2$, and all natural numbers $n$,
$n \geq 2d$, there is a polytope $Q_{n}$ that has $n$ vertices
and is contained in the Euclidean ball $B_{2}^{d}$ such that
$$
d_{H}(Q_{n},B_{2}^{d}) \leq \frac{16}{7}
(\frac{vol_{d-1}(\partial B_{2}^{d})}{vol_{d-1}( B_{2}^{d-1})})^
{\frac{2}{d-1}}
n^{-\frac{2}{d-1}}.
$$
\endproclaim
\vskip 5mm

In particular, since a $Q_n$ which satisfies the hypothesis of Lemma 3 
contains the Euclidean ball of radius $1-d_H(Q_n,B^d_2),$ it follows that
$$
d_{S}(Q_{n},B_{2}^{d})  
\leq vol_d(B^d_2)(1-d_H(Q_n,B^d_2))^d
\tag 1
$$
$$
\leq vol_{d}(B_{2}^{d})
(1-(1-\frac{16}{7}(\frac{vol_{d-1}(\partial B_{2}^{d})}{vol_{d-1}(
B_{2}^{d-1})})^
{\frac{2}{d-1}}n^{-\frac{2}{d-1}})^{d})
$$
and
$$
(1-\frac{16}{7}(\frac{vol_{d-1}(\partial B_{2}^{d})}{vol_{d-1}( B_{2}^{d-1})})^
{\frac{2}{d-1}}n^{-\frac{2}{d-1}})^{d-1}
vol_{d-1}(\partial B_{2}^{d}) \leq vol_{d-1}(\partial Q_{n}).
\tag 2
$$
We have that 
$$
vol_{d-1}(\partial B_{2}^{d})=d\ vol_{d}( B_{2}^{d})
=d \frac{\pi^{\frac{d}{2}}}{\Gamma(\frac{d}{2}+1)}
$$
$$ 
=d \sqrt{\pi}\frac{\Gamma(\frac{d-1}{2}+1)}{\Gamma(\frac{d}{2}+1)}
vol_{d-1}( B_{2}^{d-1}) 
\leq d \sqrt{\pi}\ vol_{d-1}( B_{2}^{d-1}).
$$
Since $d^{\frac{2}{d-1}} \leq 4$ and $(1-t)^{d} \geq 1-dt$ we get
from (1)
$$
d_{S}(Q_{n},B_{2}^{d}) \leq
(1-(1-\frac{64}{7}\pi n^{-\frac{2}{d-1}})^{d}) vol_{d}(B_{2}^{d}) \leq
\frac{64}{7}\pi dn^{-\frac{2}{d-1}} vol_{d}(B_{2}^{d}). 
\tag 3
$$
Similarly we get from (2) that we have for
$n \geq (\frac{128}{7}\pi d)^{\frac{d-1}{2}}$ 
$$
vol_{d-1}(\partial B_{2}^{d}) \leq 
2\ vol_{d-1}(\partial Q_{n}).  
\tag 4
$$

For the sake of completeness we include the proof of Lemma 3.
The arguments are from [BI].

\demo{Proof}
For every $n$ there is a $\theta_{n} >0$ and a set
$\{x_{1}, \dots, x_{n} \} \subset \partial B_{2}^{d}$ so that
for all $i \neq j$ we have
$$
\|x_{i}-x_{j}\| \geq \theta_{n}
$$
and so that for every $x \in \partial B_{2}^{d}$ there is $i$
such that
$$
\|x-x_{i}\| \leq \theta_{n}.
$$
We choose $Q_{n}$ to be the convex hull of $\{x_{1}, \dots, x_{n} \}$.
We have 
$$
d_{H}(Q_{n},B_{2}^{d}) \leq \frac{1}{2} \theta_{n}^{2}.
$$
If not, then there is $x \in \partial B_{2}^{d}$ such that
the Euclidean ball with radius $\frac{1}{2} \theta_{n}^{2}$
and center $x$ and $Q_{n}$ have an empty intersection. By the 
theorem of Hahn-Banach there is a hyperplane separating
$Q_{n}$ and $B_{2}^{d}(x,\frac{1}{2} \theta_{n}^{2})$. This
hyperplane cuts off a cap of height greater than 
$\frac{1}{2} \theta_{n}^{2}$. The point at the top of this cap
has a distance greater than $\theta_{n}$ from all $x_{i}, i=1,
\dots, n$. This cannot be.
\par
Now we estimate $\theta_{n}$ from above. The caps
$$
\partial B_{2}^{d} \cap H^{-}((1-\frac{1}{8}\theta_{n}^{2})x_{i}, x_{i})
\hskip 10mm
i=1, \dots, n
$$
have disjoint interiors. Therefore we get
$$
vol_{d-1}(\partial B_{2}^{d}) \geq
\sum_{i=1}^{n} vol_{d-1}
(\partial B_{2}^{d} \cap H^{-}((1-\frac{1}{8}\theta_{n}^{2})x_{i}, x_{i}))
$$
$$
\geq n (\frac{1}{2}\theta_{n} \sqrt{1-\frac{1}{16}\theta_{n}^{2}})^{d-1}
vol_{d-1}(B_{2}^{d-1}).
$$
We obtain
$$
\frac{1}{2}\theta_{n} \sqrt{1-\frac{1}{16}\theta_{n}^{2}} \leq
(\frac{1}{n}\ \frac{vol_{d-1}(\partial B_{2}^{d})}{vol_{d-1}(B_{2}^{d-1})})^
{\frac{1}{d-1}}.
$$
For $n=2d$ we get that $\theta_{n} \leq \sqrt{2}$. Indeed, just consider the set
$\{e_{1}, \dots, e_{d}, -e_{1}, \dots, -e_{d} \}$. Thus it follows
$$
\frac{\theta_{n}}{2}\sqrt{\frac{7}{8}} \leq
(\frac{1}{n}\ \frac{vol_{d-1}(\partial B_{2}^{d})}{vol_{d-1}(B_{2}^{d-1})})^
{\frac{1}{d-1}}
$$
and thus
$$
\frac{1}{2}\theta_{n}^{2} \leq \frac{16}{7}
(\frac{1}{n}\ \frac{vol_{d-1}(\partial B_{2}^{d})}{vol_{d-1}(B_{2}^{d-1})})^
{\frac{2}{d-1}}.
$$
\enddemo
\qed
\vskip 5mm

\proclaim{\smc Lemma 4}
(i) For $k=0,1,2,\dots$ and $d=1,2, \dots$ we have
$$
\int_{\Bbb{R}_{+}^{d}}(\sum_{i=1}^{d}y_{i})^{k}
\exp(-(\sum_{i=1}^{d}y_{i})^{2})dy
=
\frac{\Gamma(\frac{k+d}{2})}{2(d-1)!}.
$$
(ii)
$$
\int_{\Bbb{R}_{+}^{d}}(\sum_{i=1}^{d}y_{i}^{2})
\exp(-(\sum_{i=1}^{d}y_{i})^{2})dy
=
\frac{d^{2}}{2(d+1)!} \Gamma(\frac{d}{2}).
$$
(iii) For $i \neq j$ we have
$$
\int_{\Bbb{R}_{+}^{d}}y_{i} y_{j}
\exp(-(\sum_{i=1}^{d}y_{i})^{2})dy
=
\frac{\Gamma(\frac{d}{2})}{4(d+1)(d-1)!}. 
$$
\endproclaim
\vskip 10mm

\demo{\smc Proof}
(i)
We denote
$$
H_{t}=\{y \in \Bbb{R}_{+}^{d} | \sum_{i=1}^{d} y_{i}=t \}.
$$
Let $d_{H_{t}}$ denote the $d-1$-dimensional Lebesgue-measure
on $H_{t}$.
We have that $vol_{d-1}(H_{t})=\frac{\sqrt{d}\ t^{d-1}}{(d-1)!}$  and get
$$
\int_{\Bbb{R}_{+}^{d}}(\sum_{i=1}^{d}y_{i})^{k}
\exp(-(\sum_{i=1}^{d}y_{i})^{2})dy
=
\int_{0}^{\infty} \int_{H_{t}} \frac{t^{k}}{\sqrt{d}}\  e^{-t^{2}}d_{H_{t}}dt
$$
$$
=
\frac{1}{(d-1)!}\int_{0}^{\infty} t^{k+d-1}\ e^{-t^{2}}dt
=
\frac{1}{2(d-1)!}\int_{0}^{\infty} s^{\frac{k+d}{2}-1} e^{-s}ds
=
\frac{1}{2(d-1)!} \Gamma(\frac{k+d}{2}).
$$
(ii)
$$
\int_{\Bbb{R}_{+}^{d}}(\sum_{i=1}^{d}y_{i}^{2})
\exp(-(\sum_{i=1}^{d}y_{i})^{2})dy
=
\frac{1}{\sqrt{d}}
\int_{0}^{\infty} \int_{H_{t}}\sum_{i=1}^{d}y_{i}^{2}e^{-t^{2}}d_{H_{t}}(y)dt
$$
$$
=\frac{1}{\sqrt{d}} \int_{0}^{\infty}e^{-t^{2}}
(\frac{d}{dt} \int_{0}^{t}
\int_{H_{s}}\sum_{i=1}^{d}y_{i}^{2}d_{H_{s}}(y)ds)dt
$$
$$
=\sqrt{d} \int_{0}^{\infty}e^{-t^{2}}
(\frac{d}{dt} \int_{0}^{t}
\int_{H_{s}}y_{1}^{2}d_{H_{s}}(y)ds)dt
$$
$$
=d \int_{0}^{\infty}e^{-t^{2}}
\frac{d}{dt} 
(\int_{{\sum_{i=1}^{d}y_{i} \leq t} \atop{0 \leq y_{i}}}
y_{1}^{2}d(y))dt
$$
$$
=d \int_{0}^{\infty}e^{-t^{2}}
\frac{d}{dt} 
(\int_{0}^{t}\int_{0}^{t-y_{1}}\int_{0}^{t-y_{1}-y_{2}} \cdots
\int_{0}^{t-\sum_{i=1}^{d-1}y_{i}}y_{1}^{2}dy_{d} \cdots dy_{3}dy_{2}dy_{1})dt
$$
$$
=\frac{d}{(d-1)!} \int_{0}^{\infty}e^{-t^{2}}
\frac{d}{dt} (\int_{0}^{t} y_{1}^{2}(t-y_{1})^{d-1}dy_{1})dt
$$
$$
=\frac{d}{(d-1)!} \int_{0}^{\infty}e^{-t^{2}}
\frac{d}{dt} (t^{d+2}\int_{0}^{1} s^{2}(1-s)^{d-1}ds)dt
$$
$$
=\frac{d}{(d-1)!} \int_{0}^{\infty}e^{-t^{2}}
\frac{d}{dt} (t^{d+2}\frac{\Gamma(3) \Gamma(d)}{\Gamma(d+3)})dt
$$
$$
=\frac{2d}{(d+1)!} \int_{0}^{\infty} e^{-t^{2}}
t^{d+1}dt=
\frac{d}{(d+1)!} \int_{0}^{\infty}s^{\frac{d}{2}} e^{-s}ds=
\frac{d \Gamma(\frac{d}{2}+1)}{(d+1)!}.
$$
(iii)
$$
\int_{\Bbb{R}_{+}^{d}}y_{i} y_{j}
\ \exp(-(\sum_{i=1}^{d}y_{i})^{2})dy
=
\frac{1}{d^{2}-d} \sum_{{1 \leq k,l \leq d} \atop {k \neq l}}
\int_{\Bbb{R}_{+}^{d}}y_{k} y_{l}
\ \exp(-(\sum_{i=1}^{d}y_{i})^{2})dy
$$
$$
=\frac{1}{d^{2}-d}\int_{\Bbb{R}_{+}^{d}}
((\sum_{i=1}^{d}y_{i})^{2}-(\sum_{i=1}^{d}y_{i}^{2}))
\exp(-(\sum_{i=1}^{d}y_{i})^{2})dy
$$
By (i) and (ii) we get for the above expression
$$
\frac{1}{d^{2}-d}
(\frac{\Gamma(\frac{2+d}{2})}{2(d-1)!}-\frac{d^{2}}{2(d+1)!} 
\Gamma(\frac{d}{2}))
=
\frac{\Gamma(\frac{d}{2})}{d^{2}-d}(\frac{d}{4(d-1)!}-\frac{d^{2}}{2(d+1)!})
=
\frac{\Gamma(\frac{d}{2})}{4(d+1)(d-1)!}.
$$
\qed
\enddemo
\vskip 10mm

For the following lemma compare also [R].

\proclaim{\smc Lemma 5} Let $x_1, \dots, x_d$ be points on the Euclidean 
sphere of radius 1, $S$ the simplex $[x_1, \dots, x_d]$, and rp(S) the
radial projection of $S$, i.e. the
spherical simplex of the points $x_1, \dots, x_d$. Let $X$ be the matrix
whose columns are the vectors $x_1, \dots, x_d$. Then we have
$$
vol_{d-1}(rp(S))=
\frac{2}{\Gamma(\frac{d}{2})} |\det(X)| 
\int_{\Bbb{R}_{+}^{d}}\ \exp(-y^{t}X^{t}Xy)dy
$$
and
$$
vol_{d}([0,rp(S)])=
\frac{2}{d\Gamma(\frac{d}{2})} |\det(X)| 
\int_{\Bbb{R}_{+}^{d}} \exp(-y^{t}X^{t}Xy)dy.
$$
\endproclaim
\vskip 10mm

\demo{Proof}
We have
$$
vol_{d}(B_{2}^{d})=\frac{\pi^{\frac{d}{2}}}{\Gamma(\frac{d}{2}+1)}
$$
and
$$
\int_{\Bbb{R}^{d}} e^{-\|z\|^{2}}dz=\pi^{\frac{d}{2}}.
$$
Therefore we get
$$
vol_{d}([0,rp(S)])=
\frac{1}{\Gamma(\frac{d}{2}+1)}
\int_{\{z=t\xi| \xi \in S\ \text{and} \ t \in \Bbb{R}_{+} \}}
e^{-\|z\|^{2}}dz.
$$
Using the substitution $z=Xy$ we get that the latter expression equals
$$
\frac{1}{\Gamma(\frac{d}{2}+1)} |\det(X)| \int_{y \geq 0}
e^{-y^{t}X^{t}Xy}dy.
$$
\enddemo
\qed

\proclaim{\smc Lemma 6}
Let $x_1, \dots, x_d$ be points on the Euclidean sphere
of radius 1, $S$ the simplex $[x_1, \dots, x_d]$, and let 
$rp(S)$ be the radial projection
of the simplex $S$. Let $H$ be the hyperplane containing
the simplex $[x_1, \dots,x_d]$ and r the radius of the $d-1$-dimensional
Euclidean ball $H \cap B_{2}^{d}$. Then we have
$$
vol_{d}([0,rp(S)])-vol_{d}([0,S]) \geq
\frac{{d}^{2}}{2(d+1)}(1-\|\frac{1}{d} \sum_{i=1}^{d} x_i \|^{2})
vol_{d}([0,S])
$$
and
$$
vol_{d}([0,rp(S)])-vol_{d}([0,S]) \geq
\frac{d \sqrt{1-r^2}}{2(d+1)}(1-\|\frac{1}{d} 
\sum_{i=1}^{d} x_i \|^{2}) vol_{d-1}(S).
$$
\endproclaim
\vskip 10mm

\demo{\smc{Proof}}
By Lemma 5 we have
$$
vol_{d}([0,rp(S)])-vol_{d}([0,S])
$$
$$
=\frac{2}{d\Gamma(\frac{d}{2})} |\det(X)| 
\int_{\Bbb{R}_{+}^{d}} \exp(-y^{t}X^{t}Xy)dy
-\frac{|\det(X)|}{d!}
$$
By Lemma 5(i) with $k=0$ the last expression equals
$$
\frac{2}{d\Gamma(\frac{d}{2})} |\det(X)| 
\int_{\Bbb{R}_{+}^{d}} \exp(-y^{t}X^{t}Xy)-\exp(-(\sum_{i=1}^{d}y_{i})^{2})dy
$$
$$
=\frac{2}{d\Gamma(\frac{d}{2})} |\det(X)| 
\int_{\Bbb{R}_{+}^{d}} (\exp((\sum_{i=1}^{d}y_{i})^{2}-y^{t}X^{t}Xy)-1)
\exp(-(\sum_{i=1}^{d}y_{i})^{2})dy
$$
We use now the inequality $1+t \leq e^{t}$ and get that the above
expression is greater than or equal to
$$
\frac{2}{d\Gamma(\frac{d}{2})} |\det(X)|
\int_{\Bbb{R}_{+}^{d}}((\sum_{i=1}^{d}y_{i})^{2}-y^{t}X^{t}Xy)\exp(-(\sum_{i=1}^{
d}y_{i})^{2})dy
$$
$$
=\frac{2}{d\Gamma(\frac{d}{2})} |\det(X)|
\sum_{i,j=1}^{d}(1-<x_{i},x_{j}>)
\int_{\Bbb{R}_{+}^{d}}y_{i}y_{j}\exp(-(y^{t}y)^{2})dy
$$
Since we have $1=<x_{i},x_{i}>$ for $i=1, \dots, d$ we get by Lemma 5(iii) for
the 
above expression
$$
=\frac{2}{d\Gamma(\frac{d}{2})} |\det(X)|
\sum_{i,j=1}^{d}(1-<x_{i},x_{j}>)\frac{\Gamma(\frac{d}{2})}{4(d+1)(d-1)!} 
=\frac{1}{2(d+1)!}(d^{2}-\|\sum_{i=1}^{d}x_{i}\|^{2})|\det(X)|.
$$
\qed
\enddemo
\vskip 10mm 

\proclaim{\smc Lemma 7}
Let $A$ be a measurable subset of $B_{2}^{d}$ such that 
the center of gravity of $A$ is contained in a cap of
height $\Delta$, $\Delta \leq 1$. Then there is a cap $C$
of height $2\Delta$ so that
$$
2 \ vol_{d}(C \cap A) \geq vol_{d}(A).
$$
\endproclaim
\vskip 5mm

\proclaim{\smc Lemma 8}
Let $P_{n}$ be a simplicial polytope with vertices
$x_{1}, \dots , x_{n}$ that are elements of $\partial B_{2}^{d}$.
Let $F_{j}$, $j=1, \dots , m$ be the $d-1$-dimensional faces of
$P_{n}$, $H_{j}$ the hyperplane containing $F_{j}$, $h_{j}$
the height of the cap $B_{2}^{d} \cap H_{j}^{-}$, and $r_{j}$
the radius of $B_{2}^{d} \cap H_{j}$.
Let $\eusm N$ be the set of integers $j$ so that
$$
h_{j} \leq \frac{1}{8}
(\frac{vol_{d-1}(\partial P_{n})}{vol_{d-1}(\partial B_{2}^{d})}\
\frac{1}{4n})^{\frac{2}{d-1}}.
$$
Then we have
$$
vol_{d-1}(\bigcup_{j \in \eusm{N}} F_{j}) \leq \frac{1}{4}
vol_{d-1}(\partial P_{n}).
$$
\endproclaim
\vskip 5mm

\demo{Proof}
We put
$$
\eusm N_{i}=\{j \in \eusm N | x_{i} \in F_{j} \}
\hskip 20mm  i=1, \dots, n
$$
and
$$
\rho=\frac{1}{8}
(\frac{vol_{d-1}(\partial P_{n})}{vol_{d-1}(\partial B_{2}^{d})}\
\frac{1}{4n})^{\frac{2}{d-1}}.
$$
Since $h_{j} \leq \rho$ we have that
$\bigcup_{j \in \eusm{N}_{i}} F_{j}$ is contained in 
$B_{2}^{d}(x_{i},2\sqrt{2\rho})$. 
$\bigcup_{j \in \eusm{N}_{i}} F_{j}$ is a subset of the boundary
of the convex set $P_{n} \cap B_{2}^{d}(x_{i},2\sqrt{2\rho})$. 
Thus we get
$$
vol_{d-1}(\bigcup_{j \in \eusm{N}_{i}} F_{j}) \leq
vol_{d-1}(\partial(P_{n} \cap B_{2}^{d}(x_{i},2\sqrt{2\rho}))).
$$
Since $P_{n} \cap B_{2}^{d}(x_{i},2\sqrt{2\rho})$ is a convex 
subset of the convex set $B_{2}^{d}(x_{i},2\sqrt{2\rho})$ we get
$$
vol_{d-1}(\bigcup_{j \in \eusm{N}_{i}} F_{j}) \leq
(8\rho)^{\frac{d-1}{2}}vol_{d-1}(\partial B_{2}^{d}) \leq
\frac{1}{4n}vol_{d-1}(\partial P_{n}).
$$
Therefore we get 
$$
vol_{d-1}(\bigcup_{j \in \eusm N}F_{j})
=vol_{d-1}(\bigcup_{i=1}^{n} \bigcup_{j \in \eusm N_{i}} F_{j})
\leq \sum_{i=1}^{n} vol_{d-1}(\bigcup_{j \in \eusm N_{i}} F_{j})
\leq \frac{1}{4} vol_{d-1}(\partial P_{n}).
$$ 
\enddemo
\qed
\vskip 5mm

\proclaim{\smc Lemma 9}
Let $P_{n}$ be a simplicial polytope with vertices 
$x_{1}, \dots , x_{n}$ that are elements fo $\partial B_{2}^{d}$.
Let $F_{j}$, $j=1, \dots , m$ be the $d-1$-dimensional faces of
$P_{n}$, $H_{j}$ the hyperplane containing $F_{j}$, $h_{j}$
the height of the cap $B_{2}^{d} \cap H_{j}^{-}$, and $r_{j}$
the radius of $B_{2}^{d} \cap H_{j}$.
Assume that we have for all $j$, $j=1, \dots, m$
$$
h_{j} \leq \frac{16}{7}
(2 \frac{vol_{d-1}(\partial B_{2}^{d})}{vol_{d-1}(B_{2}^{d-1})}
)^{\frac{2}{d-1}} n^{-\frac{2}{d-1}}
$$
and assume that
$$
vol_{d-1}(\partial B_{2}^{d}) \leq 2 vol_{d-1}(\partial P_{n}).
$$
Let $\eusm M$ be the set of integers $j$ so that
$$
\|cg(F_{j})-cg(H_{j} \cap B_{2}^{d}) \|_{2} \geq 
\frac{2^{22}-1}{2^{22}}r_{j}.
$$
Then we have
$$
vol_{d-1}(\bigcup_{j \in \eusm M} F_{j}) \leq
\frac{1}{4} vol_{d-1}(\partial P_{n}).
$$
\endproclaim
\vskip 5mm

\demo{Proof}
We put
$$
\theta= \frac{16}{7}
(2 \frac{vol_{d-1}(\partial B_{2}^{d})}{vol_{d-1}(B_{2}^{d-1})}
)^{\frac{2}{d-1}} n^{-\frac{2}{d-1}}
\leq \frac{16}{7}(2d\sqrt{\pi})^{\frac{2}{d-1}}n^{-\frac{2}{d-1}}.
$$
Since $h_{j} \leq \theta$ we have for all $j$, $j=1,\dots,m$
$$
r_{j} \leq \sqrt{2\theta}.
$$
We have that $cg(F_{j})$ is contained in a cap of height $2^{-22}r_{j}$
of the $d-1$-dimensional Euclidean ball $H_{j} \cap B_{2}^{d}$. By Lemma 7
there is a subset $\tilde{F}_{j}$ of $F_{j}$ so that $\tilde{F}_{j}$ is
contained in a cap of height $2^{-21}r_{j}$ and
$$
vol_{d-1}(F_{j}) \leq 2 vol_{d-1}(\tilde{F}_{j}).
$$
Thus the diameter of $\tilde{F}_{j}$ is less than $2^{-9}r_{j} \leq
\frac{\sqrt{2\theta}}{512}$. The set of all integers $j$ such that
$x_{i} \in \tilde{F}_{j}$ is denoted by $\eusm M_{i}$. We have that
$\bigcup_{j \in \eusm M_{i}} \tilde F_{j}$ is a subset of the boundary
of the convex set $P_{n} \cap B_{2}^{d}(x_{i},2^{-9}\sqrt{2\theta})$
and has a smaller surface area than $B_{2}^{d}(x_{i},2^{-9}\sqrt{2\theta})$.
$$
vol_{d-1}(\bigcup_{j \in \eusm M_{i}} \tilde F_{j}) \leq
(\frac{\sqrt{2\theta}}{512})^{d-1}vol_{d-1}(\partial B_{2}^{d}) \leq
\frac{4d \sqrt{\pi}}{n}(\frac{\sqrt{32}}{512\sqrt{7}})^{d-1}vol_{d-1}
(\partial P_{n}).
$$
Since $d \leq 2^{d-1}$ we get that the latter expression is smaller than
$$
\frac{4 \sqrt{\pi}}{n}(\frac{\sqrt{2}}{128})^{d-1}vol_{d-1}(\partial P_{n})
\leq \frac{\sqrt{2\pi}}{32n}vol_{d-1}(\partial P_{n})
\leq \frac{1}{8n}vol_{d-1}(\partial P_{n}).
$$
Therefore we get
$$
vol_{d-1}(\bigcup_{j \in \eusm M} F_{j}) 
=vol_{d-1}(\bigcup_{i=1}^n \bigcup_{j \in \eusm M_i} F_{j})
\leq \sum_{i=1}^{n}vol_{d-1}(\bigcup_{j \in \eusm M_{i}} F_{j})
$$
$$
\leq 2\ \sum_{i=1}^{n}vol_{d-1}(\bigcup_{j \in \eusm M_{i}} \tilde F_{j})
\leq
2 \sum_{i=1}^{n}vol_{d-1}(\bigcup_{j \in \eusm M_{i}} \tilde F_{j}) \leq
\frac{1}{4}vol_{d-1}(\partial P_{n}).
$$
\enddemo
\qed
\vskip 5mm

\demo{Proof of Theorem 1}
We consider numbers of vertices $n$ such that $n \geq 
(\frac{512}{7}\pi d)^{\frac{d-1}{2}}$.
Let $P_{n}$ be a polytope with $n$ vertices so that
$vol_{d}(B_{2}^{d})-vol_{d}(P_{n})$ is minimal. Let $Q_{n}$
be a polytope with $n$ vertices so that $d_{H}(B_{2}^{d},Q_{n})$
is minimal. By Lemma 3 we have that for all j  
$$
d_{H}(B_{2}^{d},Q_{n}) \leq \frac{16}{7} 
(\frac{vol_{d-1}(\partial B_{2}^{d})}{vol_{d-1}( B_{2}^{d-1})})^
{\frac{2}{d-1}} n^{-\frac{2}{d-1}}.
$$
We consider now the convex hull of
$P_{n}$ and $Q_{n}$.
$$
P=[P_{n},Q_{n}].
$$
$P$ has at most $2n$ vertices. Its $d-1$-dimensional faces are
denoted by  $F_{j}$, $j=1, \dots , m$. $H_{j}$ is the hyperplane 
containing $F_{j}$, $h_{j}$
the height of the cap $B_{2}^{d} \cap H_{j}^{-}$, and $r_{j}$
the radius of $B_{2}^{d} \cap H_{j}$. We may assume that $P$ is
simplicial.
We have that
$$
h_{j} \leq d_{H}(B_{2}^{d},Q_{n}) \leq \frac{16}{7} 
(\frac{vol_{d-1}(\partial B_{2}^{d})}{vol_{d-1}( B_{2}^{d-1})})^
{\frac{2}{d-1}} n^{-\frac{2}{d-1}}.
$$
By the assumption on $n$ we have that
$$
h_{j} \leq \frac{1}{8}
\quad \text{and} \quad
r_{j}=\sqrt{2h_{j}-h_{j}^{2}} \leq \frac{1}{2}.
\tag 5
$$
Also we have by (4) that
$$
vol_{d-1}(\partial B_{2}^{d}) \leq 2 vol_{d-1}(\partial Q_{n})
\leq 2 vol_{d-1}(\partial P).
$$
We apply Lemma 8 and 9 to P that has at most $2n$ vertices. Thus a factor $2$
enters
the estimates. Let $\eusm L$ be the
set of integers $j$ so that
$$
\frac{1}{8}
(\frac{vol_{d-1}(\partial P_{n})}{vol_{d-1}(\partial B_{2}^{d})}\
\frac{1}{8n})^{\frac{2}{d-1}}
\leq h_{j} \leq \frac{16}{7}
(\frac{vol_{d-1}(\partial B_{2}^{d})}{vol_{d-1}(B_{2}^{d-1})}\
\frac{1}{n})^{\frac{2}{d-1}}
\tag 6
$$
and
$$
\|cg(F_{j})-cg(H_{j} \cap B_{2}^{d}) \|_{2} <
\frac{2^{22}-1}{2^{22}}r_{j}.
\tag 7
$$
We have
$$
vol_{d-1}(\bigcup_{j \in \eusm L} F_{j}) 
\geq \frac{1}{2} vol_{d-1}(\partial P).
\tag 8
$$
We apply Lemma 6
$$
vol_{d}(B_{2}^{d})-vol_{d}(P_{n}) \geq
vol_{d}(B_{2}^{d})-vol_{d}(P) \geq
\sum_{j \in \eusm L} (vol_{d}([0,rp(F_{j})])-vol_{d}([0,F_{j}]))
$$
$$
\geq \sum_{j \in \eusm L}\frac{\sqrt{1-r_{j}^{2}}}{4}
(1-\|cg(F_{j})\|_{2}^{2})vol_{d-1}(F_{j}).
$$
By (5) we have $r_{j} \leq \frac{1}{2}$ and get that the latter expression
is greater than
$$
\sum_{j \in \eusm L}\frac{1}{8}
(1-\|cg(F_{j})\|_{2}^{2})vol_{d-1}(F_{j}).
$$
We have
$$
\|cg(F_{j})\|_{2}^{2}=(1-h_{j})^{2}+\|cg(F_{j})-cg(H_{j} \cap B_{2}^{d}) 
\|_{2}^{2}.
$$
By (7) we get for $j\in \eusm L$
$$
1-\|cg(F_{j})\|_{2}^{2} \geq 1- (1-h_{j})^{2}-(\frac{2^{22}-1}{2^{22}}r_{j})^{2}
$$
$$
= 1- (1-h_{j})^{2}-(\frac{2^{22}-1}{2^{22}})^{2}(2h_{j}-h_{j}^{2})
=(2^{-21}-2^{-44})(2h_{j}-h_{j}^{2})
\geq 2^{-21}h_{j}.
$$
Therefore 
$$
vol_d(B^d_2)-vol(P)\geq
\frac{1}{2^{24}}\sum_{j \in \eusm L} h_{j}\ vol_{d-1}(F_{j}).
$$
By (6) we get that this expression is greater than
$$
\frac{1}{2^{27}}
(\frac{vol_{d-1}(\partial P)}{vol_{d-1}(\partial B_{2}^{d})}\
\frac{1}{8n})^{\frac{2}{d-1}} \sum_{j \in \eusm L} vol_{d-1}(F_{j}).
$$
By (8) this expression is greater than
$$
\frac{1}{2^{29}}
(\frac{vol_{d-1}(\partial P)}{vol_{d-1}(\partial B_{2}^{d})}\
\frac{1}{8n})^{\frac{2}{d-1}} vol_{d-1}(\partial P)
$$
$$
\geq \frac{1}{2^{36}}
vol_{d-1}(\partial B_{2}^{d})\ n^{-\frac{2}{d-1}}.
$$
\enddemo
\qed

\enddocument


\Refs
\widestnumber\key{aaaaa}

\ref
\key BI
\by E.M. Bronshtein and L.D. Ivanov
\paper The approximation of convex sets by polyhedra
\jour Siberian Mathematical Journal
\yr 1975
\vol 16
\pages 1110--1112
\endref

\ref
\key D$_{1}$
\by R. Dudley
\paper Metric entropy of some classes of sets with differentiable boundaries
\jour Journal of Approximation Theory
\yr 1974
\vol 10
\pages 227--236
\endref

\ref
\key D$_{2}$
\by R. Dudley
\paper Correction to ''Metric entropy of some classes of sets with
differentiable boundaries''
\jour Journal of Approximation Theory
\yr 1979
\vol 26
\pages 192--193
\endref

\ref
\key F--T
\by L. FejesToth
\paper \"Uber zwei Maximumsaufgaben bei Polyedern
\jour Tohoku Mathematical Journal
\yr 1940
\vol 46
\pages 79--83
\endref


\ref
\key GMR$_{1}$
\by Y. Gordon, M. Meyer, and S. Reisner
\paper Volume approximation of convex bodies by polytopes
--a constructive method
\jour Studia Mathematica
\vol 111
\yr 1994
\pages 81--95
\endref

\ref
\key GMR$_{2}$
\by  Y. Gordon, M. Meyer and S. Reisner
\paper Constructing a polytope to approximate a convex body
\jour Geometriae Dedicata
\yr 1995
\vol57
\pages 217-222
\endref

\ref
\key Gr$_{1}$
\by P.M. Gruber
\paper Volume approximation of convex bodies by
inscribed polytopes
\jour Mathematische Annalen
\yr 1988
\vol 281
\pages 292--245
\endref

\ref
\key Gr$_{2}$
\by P.M. Gruber
\paper Asymptotic estimates for best and stepwise approximation of
convex bodies II
\jour Forum Mathematicum
\vol 5
\yr 1993
\pages 521--538
\endref

\ref
\key GK
\by P.M. Gruber and P. Kenderov
\paper Approximation of convex bodies by polytopes
\jour Rend. Circolo Mat. Palermo
\yr 1982
\vol 31
\pages 195--225
\endref

\ref
\key Mac
\by A.M. Macbeath
\paper An extremal property of the hypersphere
\jour Proceedings of the Cambridge Philosophical Society
\yr 1951
\vol 47
\pages 245--247
\endref

\ref
\key M\"u
\by J.S. M\"uller
\paper Approximation of the ball by random polytopes
\jour Journal of Approximation Theory
\yr 1990
\vol 63
\pages 198--209
\endref

\ref
\key R
\by C.A. Rogers
\book Packing and Covering
\yr 1964
\publ Cambridge University Press
\endref

\ref
\key
\by
\paper
\jour
\yr
\vol
\pages
\endref

\endRefs
\bye